\documentclass[11pt]{amsart}
\usepackage{anyfontsize}
\usepackage{tikz-cd}
\usepackage[headings]{fullpage}
\usepackage{amsmath}
\usepackage{amsfonts}
\usepackage{latexsym}
\usepackage{amssymb}
\usepackage{MnSymbol}
\newtheorem{thm}{Theorem}[section]

\def\homo{\mathop{\sf Hom}}

\def\limh{\sf \overleftarrow{\mathcal{H}}}

\def\sig{\sigma, \sigma^{-1}}
  
\begin{document}
 
\title[Structurally Stable Properties of Control Systems]{Structurally Stable Properties of Control Systems}

\maketitle

\begin{center}

\author{Shiva Shankar\footnote{Chennai Mathematical Institute,
Chennai (Madras), India}}

\vspace{.25cm}
Abstract of talk dedicated to the 110th anniversary of  L.S.Pontryagin, \\
Steklov Mathematical Institute of Russian Academy of Sciences, Moscow.

\end{center}

\vspace{.5cm}

\noindent We study properties of control systems that are stable with respect to perturbations. These ideas go back to the notion of structural stability of autonomous systems due to Andronov and Pontryagin. In contrast, control systems usually admit inputs, and are therefore non-autonomous.

We study linear systems defined by partial differential or difference equations. These are systems defined over the ring $A=\mathbb{C}[\partial_1, \ldots, \partial_n]$ of partial differential operators, or the ring $B=\mathbb{C}[\sigma_1, \sigma_1^{-1}, \ldots , \sigma_n, \sigma_n^{-1}]$ of Laurent polynomials, respectively. More precisely:\\

\noindent (i) Let $P \subset A^k$ be an $A$-submodule, and suppose it is generated by $p_1, \ldots, p_\ell$. Let $p_i=(p_{i1}(\partial), \ldots ,p_{ik}(\partial))$, and let $P(\partial)$ be the $\ell \times k$ matrix whose rows are the $p_i$. Let $\mathcal{C}^\infty$ be the space of smooth functions on $\mathbb{R}^n$. Then 
\begin{equation}\begin{array}{lccc} P(\partial): &(\mathcal{C}^\infty)^k & \longrightarrow & (\mathcal{C}^\infty)^\ell \\ & f & \mapsto & P(\partial)f \end{array} \end{equation}
is an $A$-module map. The {\it distributed system} $\mathcal{B}(P)$ defined by $P$ is the kernel  
$\{f|P(\partial)f=0\}$ of the above map (it does not depend on the choice of generators for $P$ which defined the matrix $P(\partial)$, indeed $\mathcal{B}(P) \simeq \homo_A(A^k/P, ~\mathcal{C}^\infty)$). 

\vspace{.2cm}
\noindent Remark: More generally, we can replace $\mathcal{C}^\infty$ by any $A$-submodule $\mathcal{F}$ of the space $\mathcal{D}'$ of distributions on $\mathbb{R}^n$, and study the system $\homo_A(A^k/P,~\mathcal{F})$. Examples include the spaces $\mathcal{S}'$ of tempered distributions, $\mathcal{S}$ of rapidly decreasing functions, the inverse limit $\limh$ of the Sobolev spaces, etc. The answers to the questions we address depend on the choice of $\mathcal{F}$ (for instance \cite{stop}). We confine ourselves here to the space $\mathcal{C}^\infty$, and the classical result of Malgrange and Palamodov that it is an injective, cogenerating, $A$-module is of crucial importance \cite{p}. \\

\noindent (ii) Let $\mathbb{Z}^n$ be the integer lattice, and let $(\mathbb{C})^{\mathbb{Z}^n}$ be the set of all functions $w:\mathbb{Z}^n \rightarrow \mathbb{C}$. The term $\sigma_i \in B$ acts on $w$ by shift: $\sigma_i(w)(m_1, \ldots m_n)=w(m_1, \ldots, m_i+1, \ldots ,m_n)$. Composition then defines the action of a monomial, and by linearity extends to an action of $B$ on $(\mathbb{C})^{\mathbb{Z}^n}$. The {\it $n$-D system} $\mathcal{B}(P)$ defined by a submodule $P \subset B^k$ is the kernel of 
\begin{equation} P(\sig): (\mathbb{C}^k)^{\mathbb{Z}^n} \longrightarrow (\mathbb{C}^\ell)^{\mathbb{Z}^n} \end{equation}
where the $\ell$ rows of the matrix $P(\sig)$ generate $P$. Again, $\mathcal{B}(P) \simeq \homo_B(B^k/P, ~(\mathbb{C})^{\mathbb{Z}^n})$. Here it is elementary that $(\mathbb{C})^{\mathbb{Z}^n}$ is an injective cogenerating $B$-module. \\

We interpret the above formulation: let $R$ denote either the ring $A$ or $B$, and $\mathcal{F}$ denote $\mathcal{C}^\infty$ or $\mathbb{C}^{\mathbb{Z}^n}$. A system is described by some $k$ ($\mathbb{C}$-valued) attributes at various points of $\mathbb{R}^n$ or $\mathbb{Z}^n$, each such  description is an {\it evolution} of the system. A priori, perhaps the system could evolve according to any $f \in \mathcal{F}^k$, but the laws governing the system restrict the possible evolutions $f$ to a subset of $\mathcal{F}^k$. Here the set of laws governing the system is a submodule $P$ of $R^k$, and to say that $f$ must satisfy these laws is to say that it lies in the kernel of equation (1) or (2) above. \\

\noindent Example:  State space systems: the state $x: \mathbb{R} \rightarrow \mathbb{R}^\ell$ of the system evolves according to $\frac{d}{dt}  x=Xx +Uu$, where $u: \mathbb{R} \rightarrow \mathbb{R}^m$ is the input, $X$ and $U$ are $\ell \times \ell$ and $\ell \times m$ matrices with entries in $\mathbb{R}$. The possible evolutions $f=(x,u): \mathbb{R} \rightarrow \mathbb{R}^{\ell +m}$ that can occur is the kernel of 
\[ (\frac{d}{dt} I_{\ell \times \ell}-X, ~-U): (\mathcal{C}^\infty(\mathbb{R}))^{\ell+m} \longrightarrow (\mathcal{C}^\infty(\mathbb{R}))^\ell \] \hspace*{\fill}$\square$\\

\vspace{.3cm}
We wish to study perturbations of such systems, hence we need to topologise the set of all systems. In the context of (i), we need to topologize the set of all submodules $P \subset A^k$, and towards this, we need to first topologize the set $\mathcal{M}_k$ of all matrices with $k$ columns and entries from the ring $A$. We consider {\it structured} perturbations - here it means that we consider matrices with a fixed number $\ell$ of rows. Denote this subset by $\mathcal{M}_{\ell,k}$. 

Let $\mathcal{M}_{\ell,k}(d)$ be the subset of those matrices in $\mathcal{M}_{\ell,k}$ whose entries are all bounded in degree by $d$. There are ${n+d \choose n}$ monomials of degree at most $d$ in $n$ indeterminates, hence we identify $\mathcal{M}_{\ell,k}(d)$ with the $\mathbb{C}$-affine space of dimension $\ell k{n+d \choose n}$ with the Zariski topology.  For $d_1 < d_2$, $\mathcal{M}_{\ell,k}(d_1) \hookrightarrow \mathcal{M}_{\ell,k}(d_2)$ as a Zariski closed subspace, and as $d$ tends to infinity, the direct limit $\mathcal{M}_{\ell,k}$ is equipped with the direct limit topology.
A similar construction equips the set of all $\ell \times k$ matrices with entries from the Laurent polynomial ring $B$ with the Zariski topology.  
These topologies descend to submodules of $A^k$ and $B^k$, and hence to distributed and $n$-D systems respectively (details appear in \cite{sh, sr}).

We can now ask if a certain property of a distributed or $n$-D system is {\it generic} with respect to the above topology. In other words, we ask if the property holds for an open dense set of systems. \\

In this talk, I ask whether the property of being {\it controllable} is generic for  distributed systems \cite{sh}. The question whether the {\it degree of autonomy} of an $n$-D system is generic is answered  in \cite{sr}. \\

\noindent Definition\cite{w}: The distributed system $\mathcal{B}(P)$ defined by a submodule $P \subset A^k$ is controllable if for any two subsets $U_1$ and $U_2$ of $\mathbb{R}^n$ whose closures do not intersect, and any two elements $f_1$ and $f_2$ of $\mathcal{B}(P)$, there is an element $f$ in $\mathcal{B}(P)$ such that $f=f_1$ on some neighbourhood of $U_1$ and $f=f_2$ on some neighbourhood of $U_2$.

This definition generalizes the definition of a controllable state space system (of Example 1) introduced by Kalman \cite{k} in 1960, here in Moscow! 

\setcounter{section}{1}

\begin{thm}\cite{ps} The distributed system $\mathcal{B}(P)$ is controllable if and only if $A^k/P$ is torsion free.
\end{thm}

\noindent Definition The distributed system $\mathcal{B}(P)$ is {\it strictly underdetermined} if the submodule $P  \subset A^k$ can be generated by fewer than $k$ elements (i.e. $\ell < k$ in the notation of (i)). 
Otherwise, it is {\it overdetermined}.

\begin{thm} \cite{sh} A generic strictly underdetermined system is controllable, for this set of systems contains a Zariski open set. Conversely, a generic overdetermined system is uncontrollable.
\end{thm}

This follows from the following characterisation of controllability:

\begin{thm} \cite{sh} Let $P$ be a submodule of $A^k$, and let $P(\partial)$ be any $\ell \times k$ matrix whose $\ell$ rows generate $P$. Suppose the ideal $\mathfrak{i}_\ell$ of $\ell \times \ell$ minors of $P(\partial)$ is nonzero (so that $\ell \leqslant k$). Then the system $\mathcal{B}(P)$ is controllable if and only if the codimension of the variety of $\mathfrak{i}_\ell$ is greater than or equal to 2. 
\end{thm}

The paper \cite{sr} considers the important notion of the degree of autonomy of an $n$-D system, and shows that this degree is a generic property with respect to the Zariski topology. To prove this, it is first shown that while it is difficult to calculate the dimension of the variety of a specific ideal of the ring $A$, generically a variety is a complete intersection.

\end{document}